\documentclass[12pt]{article}
\input{epsf.sty}
\usepackage[dvips]{graphicx}
\usepackage{color}
\usepackage{latexsym,amsmath,amsfonts,amscd}
\usepackage{amsmath,amssymb,graphicx,amsfonts,amsthm}
\usepackage{epsfig}
\usepackage{changebar}
\usepackage{pstricks}
\usepackage{pst-plot}
\usepackage{multirow}
\usepackage{subfigure}
\usepackage{placeins}
\usepackage{amssymb}
\usepackage{bbm}
\usepackage{comment}
\usepackage{appendix}
\usepackage{makecell}
\usepackage{tabularx}
\usepackage{booktabs}
\usepackage{dashbox}
\usepackage{bm}
\usepackage{ulem}
\usepackage{stackengine}
\usepackage[colorlinks,
linkcolor=blue,
anchorcolor=red,
citecolor=green]{hyperref}

\usepackage[yyyymmdd]{datetime}

\usepackage{accents}

\allowdisplaybreaks

\newtheorem{remark}{Remark}[section]

\begin{document}

\baselineskip=1.5pc

\begin{center}
{\bf A high-order, conservative and positivity-preserving 
intersection-based remapping method between meshes with isoparametric curvilinear cells}
\end{center}


\centerline{Nuo Lei\footnote{Hua Loo-Keng Center for Mathematical Sciences,
    Academy of Mathematics and Systems Science, 
    Chinese Academy of Sciences, 
    Beijing, 100190, China. 
    E-mail: nuo\_lei@lsec.cc.ac.cn.
    Research is supported in part by 
    NSFC grant 12288201. }
\footnote{Department of Applied Mathematics, 
The Hong Kong Polytechnic University, Hung Hom, Kowloon, Hong Kong. 
Research is supported in part by the 2024 Hong Kong Scholar Program G-YZ7Q. }, 
Juan Cheng\footnote{Academy for Multidisciplinary Studies,
Capital Normal University, Beijing 100048, China. E-mail:  jcheng@cnu.edu.cn.
Research is supported in part by NSFC grant 12031001 and National Key R\&D Program
of China No. 2023YFA1009003.} and 
Chi-Wang Shu\footnote{Division of Applied Mathematics, Brown University, Providence, RI 02912, USA.
E-mail: chi-wang\_shu@brown.edu. Research is supported in part by NSF grant DMS-2309249.}}

\vspace{.15in}

\centerline{\bf Abstract}
\bigskip
This paper presents a novel intersection-based remapping method for isoparametric curvi- linear meshes within the indirect arbitrary Lagrangian-Eulerian (ALE) framework, address- ing the challenges of transferring physical quantities between high-order curved-edge meshes. Our method leverages the Weiler-Atherton clipping algorithm to compute intersections be- tween curved-edge quadrangles, enabling robust handling of arbitrary order isoparametric curves. By integrating multi-resolution weighted essentially non-oscillatory (WENO) re- construction, we achieve high-order accuracy while suppressing numerical oscillations near discontinuities. A positivity-preserving limiter is further applied to ensure physical quan- tities such as density remain non-negative without compromising conservation or accuracy. Notably, the computational cost of handling higher-order curved meshes, such as cubic or even higher-degree parametric curves, does not significantly increase compared to second- order curved meshes. This ensures that our method remains efficient and scalable, making it applicable to arbitrary high-order isoparametric curvilinear cells without compromising performance. Numerical experiments demonstrate that the proposed method achieves high- order accuracy, strict conservation (with errors approaching machine precision), essential non-oscillation and positivity-preserving.

\vspace{.05in}

\vfill 


\newpage
\baselineskip=2pc

\section{Introduction}
\label{sec1}
The arbitrary Lagrangian-Eulerian (ALE) method is a powerful framework for computational fluid dynamics, 
particularly for moving boundary and multi-material problems. 
The indirect ALE method combines the strengths of the Lagrangian and Eulerian descriptions
and comprises three key components: 
First, the Lagrangian phase, where the mesh moves in conjunction with the fluid flow 
to solve the fluid dynamics equations. Second, the rezoning (remeshing) phase, which is performed 
to enhance mesh quality. Last, the remapping phase, which involves transferring physical variables 
from the original mesh to the newly rezoned mesh.
In \cite{Cheng2007}, the authors pointed out that for a Lagrangian scheme in 2D, 
it can be at most second-order accurate if curved meshes are not used.
These limitations have spurred the development of curvilinear Lagrangian schemes and ALE methods 
\cite{Cheng2008, Dobrev2012, Anderson2018, Lieberman2019, Liu2022}, which employ higher-order curved elements 
to achieve superior geometric representation, enhancing numerical accuracy, 
and robustness under large deformations.

Nevertheless, the increased geometric complexity of curvilinear meshes introduces 
substantial challenges, particularly during the remapping phase. Robust remapping techniques 
must preserve critical properties such as conservation, high-order accuracy, and positivity. 
Overcoming these challenges is essential for advancing the capabilities of curvilinear ALE methods.
While there has been considerable research and notable progress in developing remapping methods 
for straight-line meshes, the adaptation of these methods to curvilinear meshes has been relatively under-explored.

The advection-based remapping method 
\cite{Dukowicz2000, Anderson2015, Gu2023} is a popular approach,
which solves a transport equation over a pseudo-time interval to simulate the advection of physical quantities 
from a source mesh to a target mesh.
Based on the theory of virtual element projectors, Lipnikov and his collaborators \cite{Lipnikov2019, Lipnikov2020}
have developed high-order conservative advection-based remapping method for curvilinear meshes.
This approach ensures that the remapping process preserves important physical properties
and maintains the overall accuracy of the simulation. 

The flux-based (or swept-based) remapping method \cite{Dukowicz1984, Margolin2003, Kucharik2003}
leverages the divergence theorem to convert the area integral over the swept region 
into a line integral along cell boundaries, thereby computing the fluxes of physical quantities 
across interfaces and updating the values on the target mesh.
In recent years, the authors in \cite{Lipnikov2023} have extended the swept-based remapping method 
to more general curved polygonal meshes, developing conservative and bound-preserving method.

The advection-based and flux-based approaches are highly efficient for structured meshes with consistent topology 
and small mesh motions. However, they are not suitable for scenarios 
involving large mesh movements or changes in the topological connectivity of the mesh.

The intersection-based remapping method \cite{Dukowicz1987, Kenamond2021, Shashkovy2024}, 
which computes the overlappings between old and new cells, 
offers inherent conservation properties and can handle arbitrary mesh motions. 
Due to their flexibility, these methods are particularly attractive for remeshing problems 
where the topology may differ. A wide range of intersection-based remapping methods 
have been developed, spanning from two-dimensional to three-dimensional cases
\cite{Powell2015, Lei2021, Lei2023}, from low-order to high-order methods, 
and including bound-preserving properties \cite{Barlow2016, Kucharik2003, Margolin2003}. 

These methods, primarily designed for straight-edge meshes, 
face significant challenges when applied to curved-edge meshes. 
The complexity of geometric intersections and the computational cost of solving 
high-order polynomial equations make accurate intersection calculations difficult. 
Maintaining numerical precision and consistency is crucial to ensure conservation 
and avoid distortions in the remapped fields. 
These challenges underscore the need for advanced algorithms that can efficiently 
handle curved meshes while preserving essential physical properties and avoiding numerical oscillations.

In a notable advancement, the authors in \cite{Shashkovy2024} proposed the first 
intersection-based remapping method for second-order isoparametric curvilinear meshes. 
They divided the isoparametric curvilinear cell into several sub-cells 
and identified intersections using a set-based method, which is limited to second-order curves. 
Nevertheless, this approach cannot handle higher-order curves 
and therefore cannot be extended to higher-order ALE curved mesh methods, 
which significantly limits its applicability. 

In this paper, we introduce a novel two-dimensional intersection-based remapping method 
tailored for curvilinear meshes. Our approach is capable of handling 
arbitrary high-order parameterized curves,  
thereby enabling precise representation of complex geometries. 
By applying high-order multi-resolution WENO reconstruction, 
our method achieves high-order accuracy while ensuring strict conservation and essentially non-oscillatory properties. 
Additionally, the application of a positivity-preserving limiter ensures that positivity is maintained 
without compromising high-order accuracy. This work represents a significant advancement in curvilinear 
ALE methods, offering a reliable and efficient tool for high-fidelity simulations in computational fluid dynamics.

We incorporate the Weiler-Atherton (WA) clipping algorithm, a well-established method in computer graphics, 
to address intersection-based remapping for curvilinear meshes. 
The WA algorithm efficiently computes intersection points and generates accurate boundaries 
at a relatively low computational cost, making it suitable for complex geometries. 
Its ability to handle both convex and concave polygons ensures versatility in geometric operations, 
crucial for dealing with the intricate shapes of curved-edge elements. 
This marks the first application of the WA algorithm to intersection-based remapping, 
demonstrating its effectiveness and rationality in this context. 
Our approach inherits the flexibility of intersection-based methods, 
allowing for large mesh deformations and remeshing 
without requiring identical mesh connectivity or topology. 
Compared to previous approaches \cite{Shashkovy2024}, 
our method naturally extends to higher-order curvilinear meshes 
without significant computational overhead, enhancing its applicability and efficiency.

We consider solving fluid dynamics in a cell-centered finite volume framework, 
where the cell averages of each cell are known. To achieve high-order accuracy 
in the remapping procedure, we employ the multi-resolution weighted essentially non-oscillatory (WENO) 
reconstruction approach \cite{Zhu2018}. This method utilizes hierarchical central spatial stencils 
to achieve optimal accuracy in smooth regions while effectively suppressing oscillations 
near discontinuities, with the added flexibility of arbitrary positive linear weights 
summing to one. By integrating high-order WENO reconstruction into our framework, 
we aim to achieve both high-order accuracy and essentially non-oscillatory behavior 
in handling complex situations with curvilinear meshes.

The ALE method for solving CFD problems typically involves numerous non-negative physical quantities, 
such as density and internal energy. Preserving the positivity of these quantities 
throughout the simulation is therefore crucial. Many bound-preserving or positivity-preserving 
remapping methods have been developed to address this issue, 
including those based on MOOD limiters \cite{Blanchard2016} 
or slope limiters \cite{Burton2018}. Building on our previous studies \cite{Lei2021, Lei2023}, 
we extend a positivity-preserving modification \cite{Zhang2010} to the high-order 
reconstructed polynomials on curvilinear meshes. This ensures that the remapped results remain positive 
while maintaining both conservation property and high-order accuracy.

This article begins by introducing the setup of curvilinear meshes 
and outlining the overall framework of our remapping method. 
Subsequently, from Section \ref{sec3} to Section \ref{sec6}, 
we systematically elaborate on the key components of our remapping method, 
including high-order WENO reconstruction, the Weiler-Atherton clipping algorithm, 
two types of high-order numerical integrals on curved intersections, 
and the positivity-preserving modification. 
Finally, in Section \ref{sec7}, we present numerical experiments that 
demonstrate the properties of high-order accuracy, essentially non-oscillatory behavior, 
conservation and positivity-preserving  of our remapping approach. 
This comprehensive development and validation process highlights 
the robustness of our proposed method for applications involving curvilinear meshes.

\bigskip

\section{Basic concepts}
\label{sec2}
\setcounter{table}{0}
\setcounter{figure}{0}
Suppose the two-dimensional computational domain $\Omega$ is a connected domain decomposed into
isoparametric elements $\mathcal{M} = \{I_i \}_{i=1}^N$ called ``curvilinear quadrangles''.
Each second-order curvilinear quadrangle has four vertices and four control points on the edge.
For example, consider a curved-edge with two vertex coordinates 
\( P_1(x_{P_1}, y_{P_1}) \) and \( P_2(x_{P_2}, y_{P_2}) \), 
and a control point \( P_5(x_{P_5}, y_{P_5}) \). 
We construct a quadratic curve as:
\[
\begin{cases}
x(t) = x_{P_1}(1-2t)(1-t) + 4x_{P_5}t(1-t) + x_{P_2}(2t-1)t, \\
y(t) = y_{P_1}(1-2t)(1-t) + 4y_{P_5}t(1-t) + y_{P_2}(2t-1)t,
\end{cases}
\]
where \( t \in [0, 1] \). This curve is designed such that it passes through 
\( P_1 \), \( P_5 \), and \( P_2 \) at \( t = 0 \), \( t = 0.5 \), and \( t = 1 \), respectively. 
Similarly, each third-order curvilinear quadrangle has four vertices and eight control points on the edge.
Figure \ref{quad} shows some examples for the second-order and third-order curvilinear quadrangles.

For the curvilinear mesh $\mathcal{M}$, 
we have some basic requirements:
\begin{itemize}
\item There are neither overlaps nor gaps between any two neighboring cells,
$$
\bigcup_{i=1}^N I_i = \Omega, \quad 
| I_i \bigcap I_{j} | = 0 \; \forall \;  i\neq j.
$$
\item All of the cells $\{ I_i \}_{i=1}^N$ 
should not be self-intersecting.
\item All of the cells $\{ I_i \}_{i=1}^N$  
do not have holes.
\end{itemize}
Obviously, these requirements are natural and we also require the new rezoned mesh 
$\widetilde{\mathcal{M} } =  \{ \tilde{I}_{\tilde{i}} \}_{\tilde{i}=1}^{\tilde{N}}$
satisfies them.
We denote the intersection of the old cell $I_i$ and the new cell $\tilde{I}_{\tilde{i}}$
as $\hat{I}_{i, \tilde{i}} := I_i \bigcap \tilde{I}_{\tilde{i}}$.

In this work, we focus on remapping cell averages from $\mathcal{M}$ to $\widetilde{\mathcal{M} }$.
Suppose we have cell averages $\bar{u}_i$ in each cell $I_i$, then we need to calculate the new cell average 
$\bar{ \tilde{u} }_{\tilde{i}}$ in $\tilde{I}_{\tilde{i}}$ for $\tilde{i} = 1, \cdots , \tilde{N}$.
Refer to our previous work \cite{Lei2021}, the flowchart of our curved remapping method follows as, 
\begin{enumerate}
\item {\bf High-order reconstruction: }
Reconstruct high-order polynomial $u_i(x, y)$ in each old cell $I_i$ with the multi-resolution WENO idea.
\item {\bf Weiler-Atherton clipping: }
Find the intersection $\hat{I}_{i, \tilde{i}}$ between each old 
and new isoparametric cells.
\item {\bf The Positivity-preserving modification: }
Modify $u_i(x, y)$ by the positivity-preserving limiter to maintain positivity for the physically positive
variables, denote the resulting modified polynomial as $\hat{u}_i(x, y)$.
\item {\bf Numerical integration: }
Calculate the integration for the modified polynomials over the intersections 
$\int_{\hat{I}_{i, \tilde{i}} } \hat{u}_i(x, y) dxdy$.
\end{enumerate}
Finally, we can calculate the new cell averages as 
$$
\bar{\tilde{u}}_{\tilde{i}} 
:= \frac{1}{ |\tilde{I}_{\tilde{i}} | } \sum_{i=1}^N 
\int_{\hat{I}_{i, \tilde{i}} } \hat{u}_i(x, y) dxdy .
$$
Based on clipping the intersection exactly and high-order quadrature rule, 
we can prove that this curved remapping method is conservative, see \cite{Lei2021} for the detailed proof.
Next, we will present a step-by-step explanation of our remapping method.

\begin{figure}[htbp!]
	\centering
	\includegraphics[width=7.5cm]{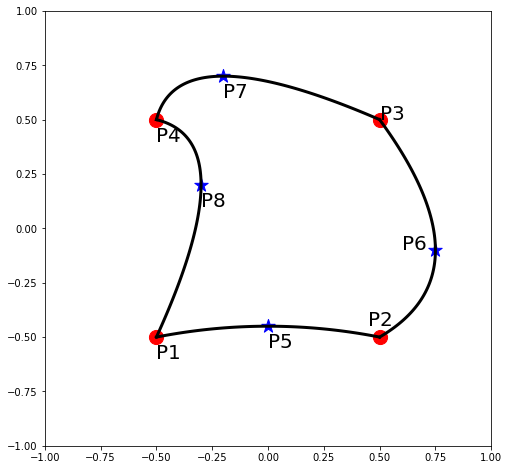}
	\includegraphics[width=7.5cm]{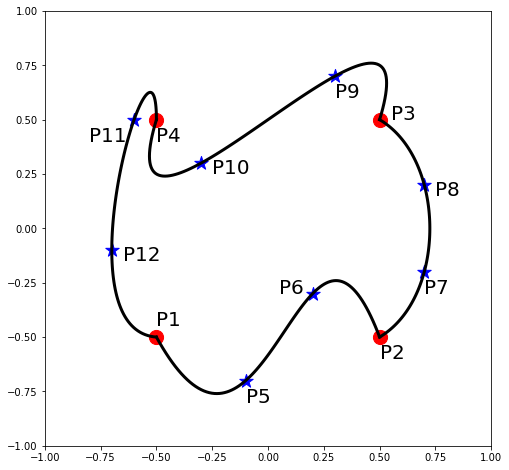}
	\caption{The second-order and third-order isoparametric cells. 
	Red markers are vertices and blue markers are control points.}
	\label{quad}
\end{figure}

\bigskip

\section{High-order reconstruction}
\label{sec3}
\setcounter{table}{0}
\setcounter{figure}{0}
\setcounter{equation}{0}

On the original curvilinear mesh \( \mathcal{M} \), 
we employ the multi-resolution WENO approach \cite{Zhu2018} to 
reconstruct high-order and conservative polynomials \( \{ u_i(x,y) \}_{i=1}^N \).
For illustration, we focus on the third-order WENO reconstruction. 
This method can be naturally extended to higher orders by increasing the number of stencils 
and polynomial degrees while maintaining the essential properties of the reconstruction.
\begin{enumerate}
\item Select a small stencil $S_0 = \{ I_i \}$ and a big stencil 
$S_1 = \{I_j | I_j\bigcap I_i \neq \varnothing\}$.
Here, we consider curved quadrilateral cells. 
Thus, the big central stencil $S_1$ includes the cell $I_i$ itself, four cells sharing edges with $I_i$,
and four cells sharing vertices with $I_i$.
It can also be seen that $S_0$ is contained within $S_1$, $S_0 \subset S_1$.
\item Reconstruct a zeroth-degree polynomial $q_0(x,y) = \bar{u}_i$ on $S_0$
and reconstruct quadratic polynomial $q_1(x,y) \in \mathbb{P}^2(I_i)$ on  $S_1$ such that 
\begin{equation}\label{leastsquare}
\begin{aligned}
q_1(x,y) &= \underset{\tilde{q}_1\in \mathbb{P}^2}{\arg\min} \sum_{I \in S_1} 
\left| \iint_{I} \tilde{q}_1(x,y) dxdy - \bar{u}_I |I| \right|^2 , \\
&\text{s.t.}\; \iint_{I_i} \tilde{q}_1(x,y) dxdy = \bar{u}_{i} |I_i| .
\end{aligned}
\end{equation}
where $|I|, |I_i|$ are the area of the relative cells.
\item Take $p_0(x,y) = q_0(x,y)$ and 
$p_1(x,y) = \frac{1}{\gamma_1}q_1(x,y) - \frac{\gamma_0}{\gamma_1}p_1(x,y)$,
with linear weights $\gamma_0 + \gamma_1 = 1$ and $\gamma_0, \gamma_1 > 0$.
Following \cite{Zhu2018}, we choose $\gamma_0 = \frac{1}{101}$ and $\gamma_1 = \frac{100}{101}$.
\item Define the smoothness indicators as 
\begin{equation}\label{beta}
\beta_1 := \sum_{l_1+l_2 = 1, 2}\iint_{I_i} |I_i|^{l_1 + l_2 - 1} 
\left(\frac{\partial^{l_1 + l_2}}{\partial x^{l_1} \partial y^{l_2}} p_1(x,y)\right)^2 dxdy,\quad 
l_1, l_2 = 0, 1, 2,
\end{equation}
and $\beta_0 = \min\limits_{I_j \in S_1 \backslash \{I_i\} }\{ (\bar{u}_i - \bar{u}_j)^2 \}$, 
where the designed smoothness indicators $\beta_0, \beta_1$ 
can avoid the numerical oscillation near the discontinuity.
\item Define the nonlinear weights as 
$$
\omega_0 = \frac{\bar{\omega}_0 }{\bar{\omega}_0 + \bar{\omega}_1}, \; 
\omega_1 = \frac{\bar{\omega}_1 }{\bar{\omega}_0 + \bar{\omega}_1}, \; 
\text{with}\; \bar{\omega}_l = \gamma_l\left(1 + \frac{\tau}{\beta_l + \epsilon} \right), \;
l = 0, 1, 
$$
where $\tau = \frac{|\beta_0 - \beta_1|^2}{4}$ and $\epsilon = 10^{-4}$. 
\item Finally, the reconstructed polynomial is 
$$ 
u_i(x,y) := \omega_0 p_0(x,y) + \omega_1 p_1(x, y).
$$
\end{enumerate}
Notice that, in \eqref{leastsquare} and \eqref{beta}, 
the integrals of polynomials over curved cells are needed. 
We use Green's theorem to transform the integral over the curved cells 
into a line integral, 
and then apply one-dimensional high-order integration formulas for its computation.
For example, $\partial I_i $ has four segments $\{L_l\}_{l=1}^4$ 
and each segment \( L_l \) is parameterized by the expressions 
\( x = x_l(t),  y = y_l(t) \) for $t\in [0, 1]$. 
Then, we have, 
\begin{equation}\label{Integral_Green1}
\begin{aligned}
\iint_{I_{i}} f(x,y) dxdy 
&= \int_{\partial I_{i}} f_1(x,y)dx + f_2(x,y) dy \\
&= \sum_{l=1}^{4} \int_{0}^{1} \left(f_1(x_l(t), y_l(t)) x_l'(t) + f_2(x_l(t), y_l(t)) y_l'(t)\right)  dt
\end{aligned}
\end{equation}
where $f(x,y)$ is the basis polynomial integrand and 
polynomials $f_1(x,y), f_2(x,y)$ satisfy $\frac{\partial f_2}{\partial x} - \frac{\partial f_1}{\partial y} = f$.

In the multi-resolution WENO approach, 
the positive linear weights are flexible, and fewer high-order reconstructions are required, 
making the method significantly more efficient than the classical WENO approach. 
In the procedure, 
we  reconstruct a high-order polynomial and a low-order polynomial 
on the central nested stencils, 
with nonlinear weights assigned based on their smoothness.
This reconstruction ensures that the resulting polynomial retains high-order accuracy in smooth regions 
while placing greater emphasis on low-order polynomials near discontinuities, 
effectively suppressing numerical oscillations and improving algorithm stability.

\bigskip

\section{The Weiler-Atherton clipping algorithm}
\label{sec4}
\setcounter{table}{0}
\setcounter{figure}{0}
\setcounter{equation}{0}
The Weiler-Atherton (WA) clipping algorithm \cite{WA} is a robust computational geometry technique 
for clipping polygons. Unlike the simpler methods such as the Sutherland-Hodgman algorithm \cite{Sutherland1974}, 
which are limited to convex polygons and straight-edge meshes, 
the WA algorithm supports both convex and concave polygons, 
making it highly versatile for isoparametric curved meshes. 
It works by identifying intersection points between the subject and clipping polygons, 
then traversing their boundaries to construct the clipped result to ensure accurate traversal order. 
The flowchart of the WA clipping algorithm follows as,

\begin{enumerate}
\item {\bf Intersection points detection and labeling:}
First, we calculate all intersection points between the curved edges. 
Each intersection point is then labeled as either ``entering'' or ``exiting'' the clipping polygons, 
depending on the relative orientation of the two polygons. 
To achieve this, we use nonlinear root-finding algorithms, 
such as the Newton method, 
to solve the parametric equations 
between two curved edges and locate the intersection points.
In our implementation, the error tolerance for the nonlinear root-finding algorithm 
is set to \( 10^{-12} \). 

\item
{\bf Path construction: }
Once all intersection points are determined and labeled, 
the algorithm constructs the resultant polygon(s). 
This is done by traversing the subject polygon and clipping polygon alternately, 
starting from an ``entering'' intersection point and following the edges of the polygons 
until an ``exiting'' point is reached. 
At each intersection, the algorithm switches between the subject and clipping polygons. 
This ensures that only the portions of the subject polygon inside the clipping polygon are retained. 
The traversal continues until a closed path is formed. 
For complex cases, such as multiple disjoint intersections, 
the algorithm handles each path independently, 
generating multiple output polygons as needed.
\end{enumerate}

A key strength of the WA algorithm is its ability to handle complex scenarios, 
such as multiple intersections or disjoint intersections. 
Notably, it is also applicable to curvilinear polygons, making it suitable for finding intersections 
in curved meshes.  
Despite these constraints, its flexibility and accuracy make it a valuable tool in computer graphics, 
CAD systems, and GIS for complex geometric operations.

\begin{figure}
	\centering
	\includegraphics[width=12cm]{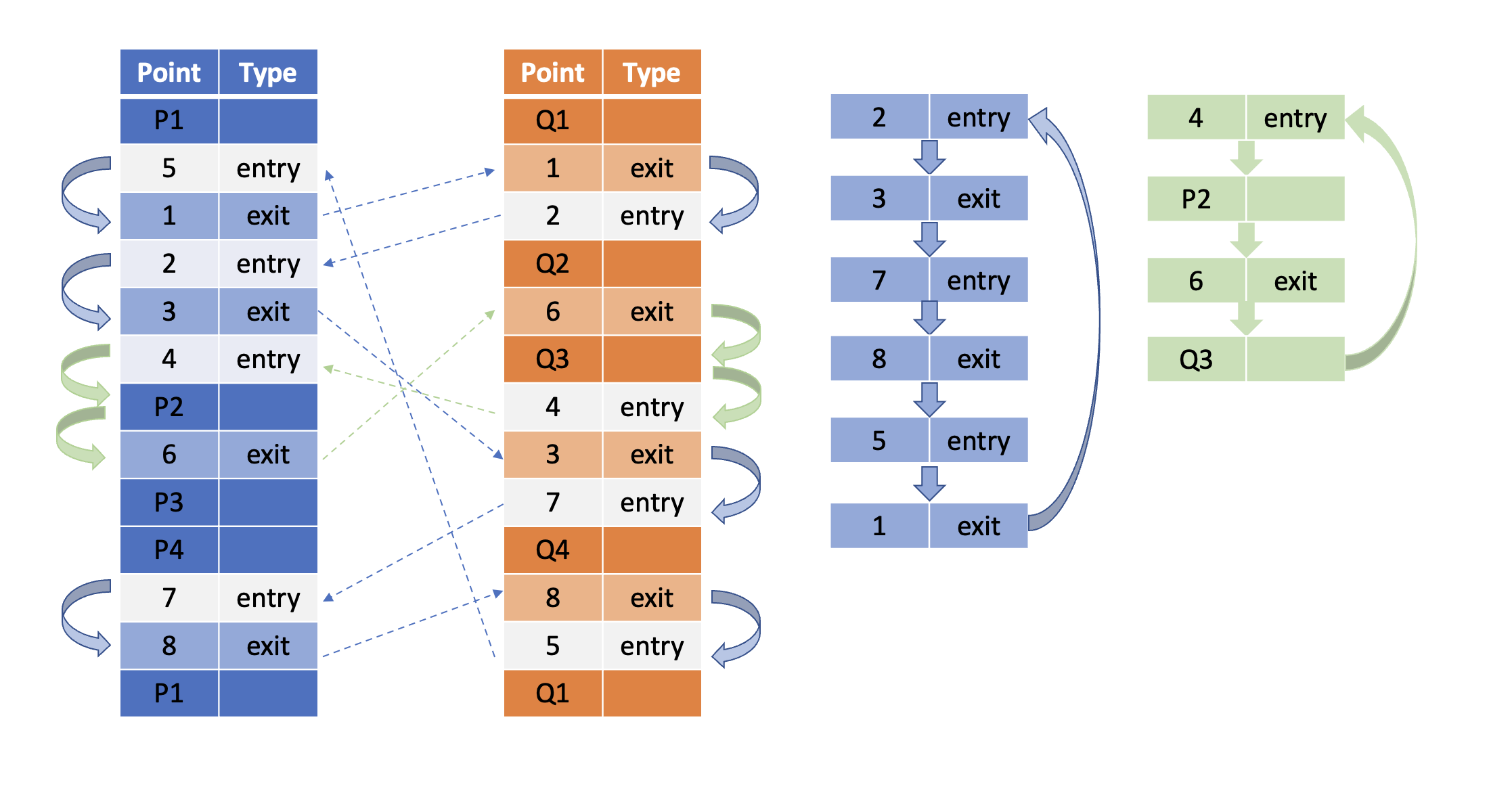}
	\caption{Linked list example.}
	\label{WA_example}
\end{figure}

Following \cite{Shashkovy2024}, we consider two cells ``Quad P'' and ``Quad Q''
in Table \ref{QuadExample},
and each of the curves is a quadratic curve.
Specifically, we consider using this example to illustrate our WA algorithm. 
\begin{enumerate}
\item [Step 1:] Define the contours for the two quadrilaterals in a counterclockwise order
$P_1 \rightarrow P_2 \rightarrow P_3 \rightarrow P_4 \rightarrow P_1, \;
Q_1 \rightarrow Q_2 \rightarrow Q_3 \rightarrow Q_4 \rightarrow Q_1$.
Based on the vertices and control points, we derive the parametric expressions 
$\{x^P_l(t), y^P_l(t) \}_{l=1}^4$ for ``Quad P'' with $t\in [0, 1]$
and $\{x^Q_l(s), y^Q_l(s) \}_{l=1}^4$ for ``Quad Q'' with $s\in [0, 1]$. 
\item [Step 2:] Compute the intersections of the curves pairwise, 
using the Newton method to solve the nonlinear equations 
$$
\left\{ \begin{matrix}
x^P_{l_1}(t) = x^Q_{l_2}(s) \\
y^P_{l_1}(t) = y^Q_{l_2}(s)
\end{matrix}\quad \text{with}\; t, s\in[0,1] \quad \forall 1\leq l_1, l_2\leq 4, \right.
$$
to obtain the coordinates of intersection points. 
Then, label the entry or exit of the intersection according to the connecting order 
of the contours, as shown in Table \ref{CoordinateExample}.
\item [Step 3:] Place the intersection points in sequential order into our linked lists, 
resulting in two new linked lists, as shown on the left side of Figure \ref{WA_example}. 
\item [Step 4:] Start from the entry point of the subject polygon 
and traverse along the linked list. If an exit point is encountered, 
switch to the clipping polygon. Continue traversing along the linked list of the clipping polygon, 
and if an entry point is encountered, switch back to the subject polygon. 
This process of switching between the two linked lists continues until a closed loop is formed. 
Since there may be multiple intersections, 
repeat this path selection process until all intersection points have been utilized.
See the right side of Figure \ref{WA_example} for the two closed loops.
\end{enumerate}

\begin{figure}
	\centering
	\includegraphics[width=7.5cm]{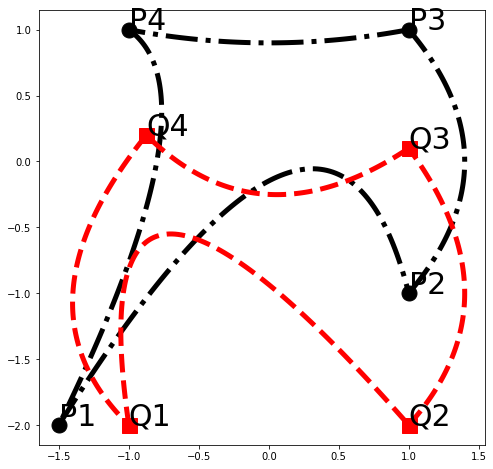}
	\includegraphics[width=7.5cm]{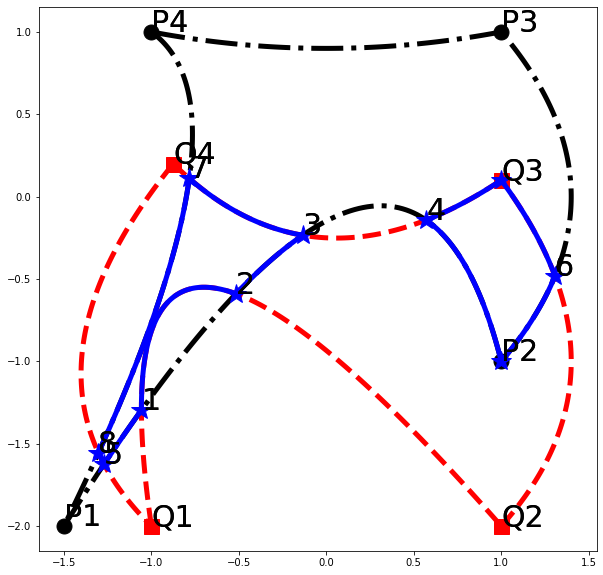}
	\caption{Clipping example with the WA algorithm.}
	\label{WA_test1}
\end{figure}

\begin{figure}
	\centering
	\includegraphics[width=5cm]{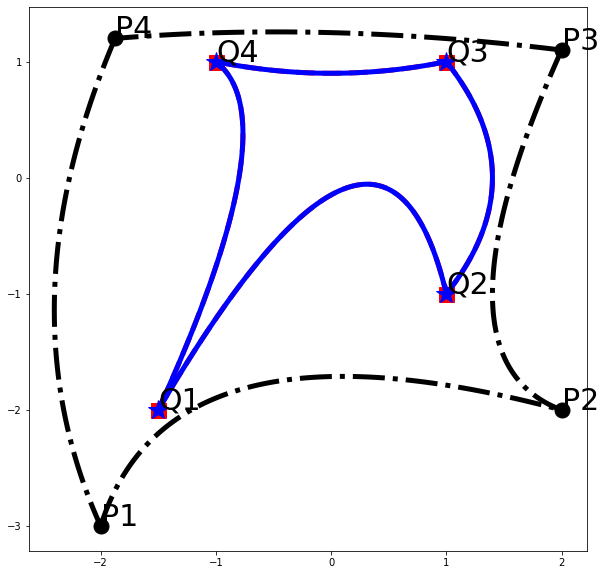}
	\includegraphics[width=5cm]{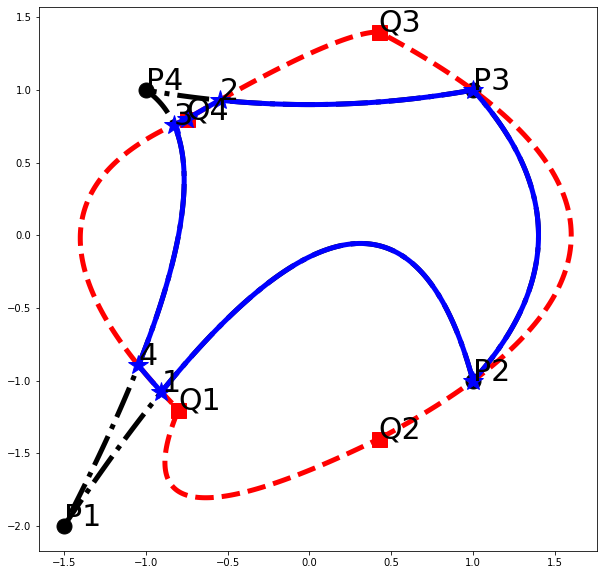}
	\includegraphics[width=5cm]{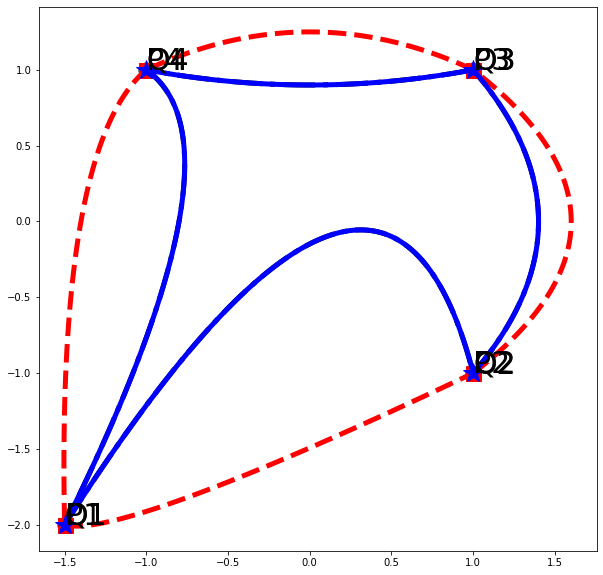}
	\includegraphics[width=5cm]{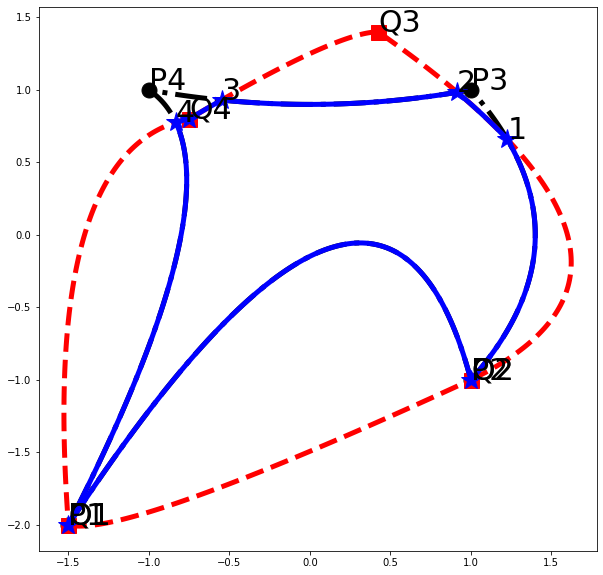}
	\includegraphics[width=5cm]{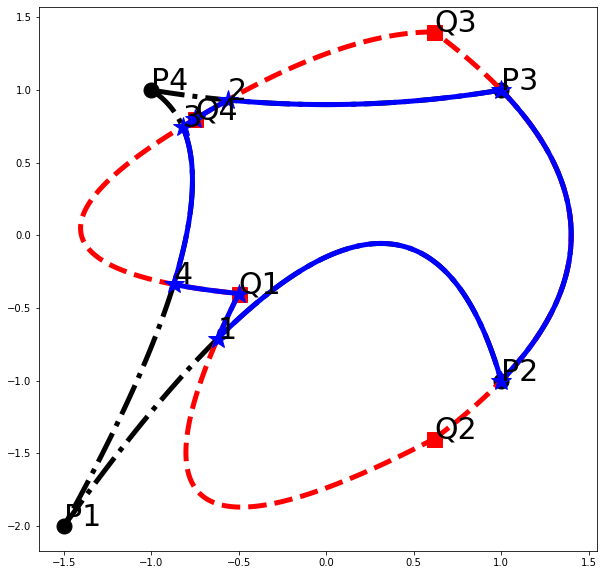}
	\includegraphics[width=5cm]{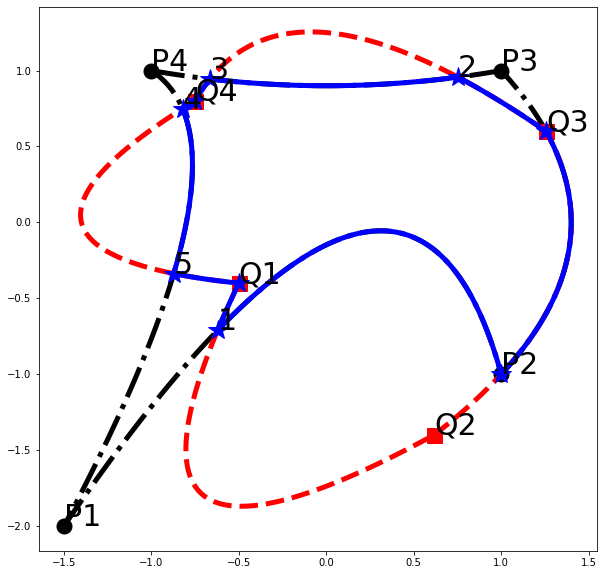}
	\caption{Some degenerate cases with the WA clipping algorithm.}
	\label{WA_test2}
\end{figure}

In Figure \ref{WA_test1}, we mark the two cells in black and red lines, 
and mark the intersections in blue solid lines.
Moreover, in practical applications, careful handling of various degenerate cases is necessary, 
such as point-to-point, point-to-edge, and edge-to-edge coincidences. 
In Figure \ref{WA_test2}, we present several examples of the WA algorithm under extreme conditions, 
demonstrating its robustness.

\begin{figure}
	\centering
	\includegraphics[width=7.5cm]{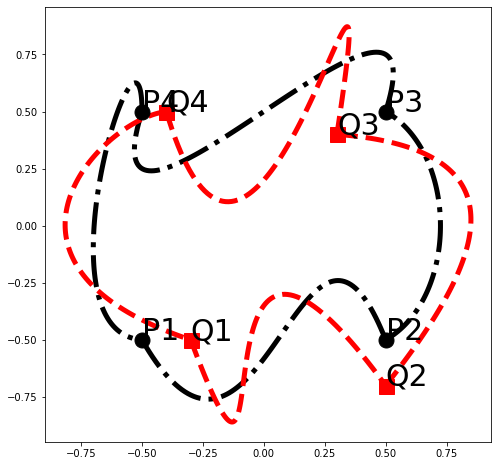}
	\includegraphics[width=7.5cm]{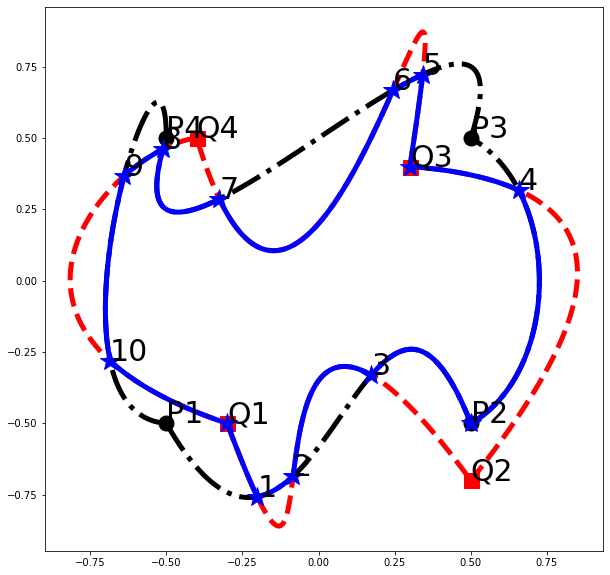}
	\caption{Clipping example with the WA algorithm for third-order isoparametric cells.}
	\label{WA_test_3rd}
\end{figure}

\begin{remark}
It can be seen from the above flowchart that the WA clipping algorithm 
can be easily extended to higher-order curved meshes 
or other types of curved grids without a significant increase in computational cost. 
For example, we demonstrate the computational results for third-order curved meshes 
in Figure \ref{WA_test_3rd}. 
\end{remark}

\begin{table}[ht]
\centering
\caption{Example quadratic cells.}
\begin{tabular}{c|c|c||c|c|c}
\toprule
\multicolumn{3}{c||}{ Quad P } & \multicolumn{3}{c}{Quad Q} \\ \hline
\quad & $x$ & $y$ & \quad & $x$ & $y$ \\ \hline
$P_1$ & -1.5 & -2 & $Q_1$ & -1 & -2\\ \hline
$P_2$ & 1 & -1 & 	$Q_2$ & 1 & -2\\ \hline
$P_3$ & 1 & 1 & 	$Q_3$ & 1 & 0.1\\ \hline
$P_4$ & -1 & 1 & 	$Q_4$ & -0.875 & 0.2\\ \hline
$P_5$ & 0.1 & -0.1 & $Q_5$ & -0.7 & -0.55\\ \hline
$P_6$ & 1.4 & 0 & $Q_6$ & 1.4 & -1\\ \hline
$P_7$ & 0 & 0.9 & $Q_7$ & 0 & -0.25\\ \hline
$P_8$ & -0.8 & 0 & $Q_8$ & -1.4 & -1\\ 
\bottomrule
\end{tabular}
\label{QuadExample}
\end{table}

\begin{table}[ht]
\centering
\caption{Coordinates of the intersection points of Quad A and Quad B in Table \ref{QuadExample}. }
\begin{tabular}{c|c|c|c||c|c|c|c}
\toprule
\multicolumn{8}{c}{ Intersection points } \\ \hline
\quad & $x$ & $y$ & type & \quad & $x$ & $y$ & type \\ \hline
`1' & -1.05713663 & -1.29598616 & exit & `2' & -0.51624632 & -0.58935822 & entry \\
`3' & -0.13487741 & -0.23429411 & exit & `4'& 0.57041875 & -0.14211494 & entry\\
`5' & -1.27097304 & -1.62405364 & entry& `6' & 1.30701824 & -0.48213526 & exit \\
`7' & -0.78370680 & 0.11026931 & entry & `8' & -1.30299176 & -1.55541394 & exit \\
\bottomrule
\end{tabular}
\label{CoordinateExample}
\end{table}

\bigskip

\section{Numerical integration over curvilinear intersection}
\label{sec5}
\setcounter{table}{0}
\setcounter{figure}{0}
\setcounter{equation}{0}
After clipping between curved cells, we next address the computation of high-order polynomial integrals 
over the resulting intersections. In \cite{Lei2021}, where the intersections are standard polygons, 
they are subdivided into multiple triangles by connecting vertices, 
and a high-order quadrature rule is applied to evaluate the integrals. 
However, this method is not directly applicable to curved-edge polygons. 
When connecting two vertices with a straight line, it may intersect other curved edges, 
resulting in an incorrect partitioning. 
To accurately compute the integral of a polynomial over the intersection of curved-edge polygons, 
in this section, we will propose two approaches.

First, we transform the two-dimensional integral into a one-dimensional contour integral 
along the curved boundaries, which has been mentioned in \eqref{Integral_Green1}. 
Specifically, assume the edge of the intersection 
$\partial \hat{I}_{i, \tilde{i}}$ has \( S \) segments $\{L_l\}_{l=1}^S $ 
with expressions 
\( x = x_l(t),  y = y_l(t) \) for $t\in [t_l^{\rm sta}, t_l^{\rm end}]$ and $l=1,\cdots, S$. 
Notice that, since each edge of the intersection lies either on $\partial I_i$ or $\partial \tilde{I}_{\tilde{i}}$, 
the expressions \((x_l(t), y_l(t))\) are known, and the parameters \(t_l^{\rm sta}, t_l^{\rm end}\) 
are already computed during the WA clipping algorithm. 
This ensures that all necessary geometric information is readily available for further calculations.
Then, we use Green's theorem such that, 
\begin{equation}\label{Integral_Green2}
\begin{aligned}
\iint_{\hat{I}_{i, \tilde{i}}} f(x,y) dxdy 
= \sum_{l=1}^{S} \int_{t_l^{\rm sta} }^{t_l^{\rm end}}\left( f_1(x_l(t), y_l(t)) x_l'(t) + f_2(x_l(t), y_l(t)) y_l'(t)\right)  dt
\end{aligned}
\end{equation}
where $f(x,y)$ is the polynomial integrand and 
polynomials $f_1(x,y), f_2(x,y)$ satisfy $\frac{\partial f_2}{\partial x} - \frac{\partial f_1}{\partial y} = f$.
Setting \( f = 1 \), 
we could compute the area of the curved-edge intersection using the following formula: 
\begin{align*}
\left| \hat{I}_{i, \tilde{i}} \right|
= \iint_{ \hat{I}_{i, \tilde{i}} } 1 dxdy 
= \frac{1}{2} \int_{t_l^{\rm sta}}^{t_l^{\rm end}} ( x_l(t)y'_l(t) - y_l(t)x'_l(t)) dt .
\end{align*}

The second approach involves partitioning the curved-edge polygons into multiple small curved-edge triangles 
and applying high-order quadrature rules for integration. 
As previously mentioned, when the intersection involves curved-edge polygons, 
directly decomposing the intersection 
into small triangles to locate all integration points becomes nontrivial.

The ear-clipping algorithm \cite{Meisters1975} is a classical method for triangulating polygons 
and is applicable to arbitrary 2D polygons, including both convex and concave shapes. 
This method works by iteratively identifying and removing ``ears" 
(triangles formed by three consecutive vertices that do not contain any other vertices of the polygon), 
ultimately decomposing the polygon into a set of non-overlapping triangles. 

Therefore, we first insert a set of points \( \{ v^n_l \}_{l=1}^S, n \in \mathbb{N} \) 
along the intersection edges \( \partial \hat{I}_{i, \tilde{i}} \)
and employ the ear-clipping algorithm  
to iteratively extract small curved-edge triangles by connecting vertices 
and inserting additional points $ \{ v^n_l \}_{l=1}^S $, 
thereby ensuring a complete and accurate partitioning of the intersection.
If a self-intersecting curved triangle appears, 
we will continue to increase points until no self-intersecting curved triangle is produced.

Next, we only need to compute the integral of the polynomial over each curved-edge triangle $ \mathbb{T} $.
Note that the curved edge of $\mathbb{T}$ must correspond to certain segments \( \{L_l \}_{l=1}^S \). 
So, we can accurately identify the control points on each edge of \( \mathbb{T} \). 
Using the three vertices and the three middle control points of \( \mathbb{T} \), 
we establish a bijection \( F_{\mathbb{T}}: \mathbb{T}_0 \rightarrow \mathbb{T} \), 
where $ {\mathbb{T}}_0 = \{(0,0), (0,1), (1,0) \} $ is the reference triangle. 
Assume we have the following quadrature rule on $\mathbb{T}_0$, 
$$
\int_{\mathbb{T}_0} f(\xi, \eta) d\xi d\eta = \sum_{m = 1}^{M} \omega_m f(\xi_m, \eta_m),
$$
where $\{ \omega_m \}_{m = 1}^M$ are the quadrature weights 
and $\{ (\xi_m, \eta_m) \}_{m = 1}^M$ are the quadrature points in $\mathbb{T}_0$.
Note that we need to require the quadrature weights be positive $\omega_m > 0$, 
which is very important for us to preserve positivity in the next section.
Then the integral over the curved-edge triangle \( \mathbb{T} \) can be expressed as  
\begin{equation}\label{Integral_par}
\begin{aligned}
\iint_{\mathbb{T} } f(x, y)  dxdy 
&= \iint_{\mathbb{T}_0} f(F_{\mathbb{T}}(\xi, \eta)) \, J_{F_{\mathbb{T}}}(\xi, \eta) \, d\xi d\eta \\
&= \sum_{m = 1}^M \omega_m f(F_{\mathbb{T}}(\xi_m, \eta_m)) \, J_{F_{\mathbb{T}}}(\xi_m, \eta_m) 
\end{aligned}
\end{equation}  
where \( J_{F_{\mathbb{T}}}(\xi, \eta) \) is the Jacobian determinant of the mapping $F_{\mathbb{T}}$.
In order to integrate a polynomial $f(x,y)$ of degree $k$ on the curved cell $\mathbb{T}$ exactly,
we would need the quadrature rule in \eqref{Integral_par} on the reference cell $\mathbb{T}_0$ to
be exact for all polynomials of degree up to $2k+2$.

\begin{remark}
The process of partitioning a polygonal into triangles 
is referred as ``triangulation'' in computer graphics. 
Evidently, the more points inserted along the curved edges, 
the more accurately the intersection of the curved-edge polygon 
can be triangulated, even though some resulting triangles may be curved-edge triangles. 
However, inserting too many points may reduce computational efficiency. 
To address this, curvature-adaptive or arc-length-adaptive inserting methods 
can be employed to insert points, making triangulation more efficient. 
Additionally, inserting points within the intersection can produce higher-quality triangles, 
as seen in methods like the Constrained Delaunay Triangulation (CDT) \cite{Chew1989}. 
That said, how to efficiently triangulate is beyond the scope of this paper. 
Instead, we focus on ensuring that sufficient tools are available 
to enable the triangulation of curved-edge intersection.
\end{remark}

\begin{figure}
	\centering
	\includegraphics[width=7.5cm]{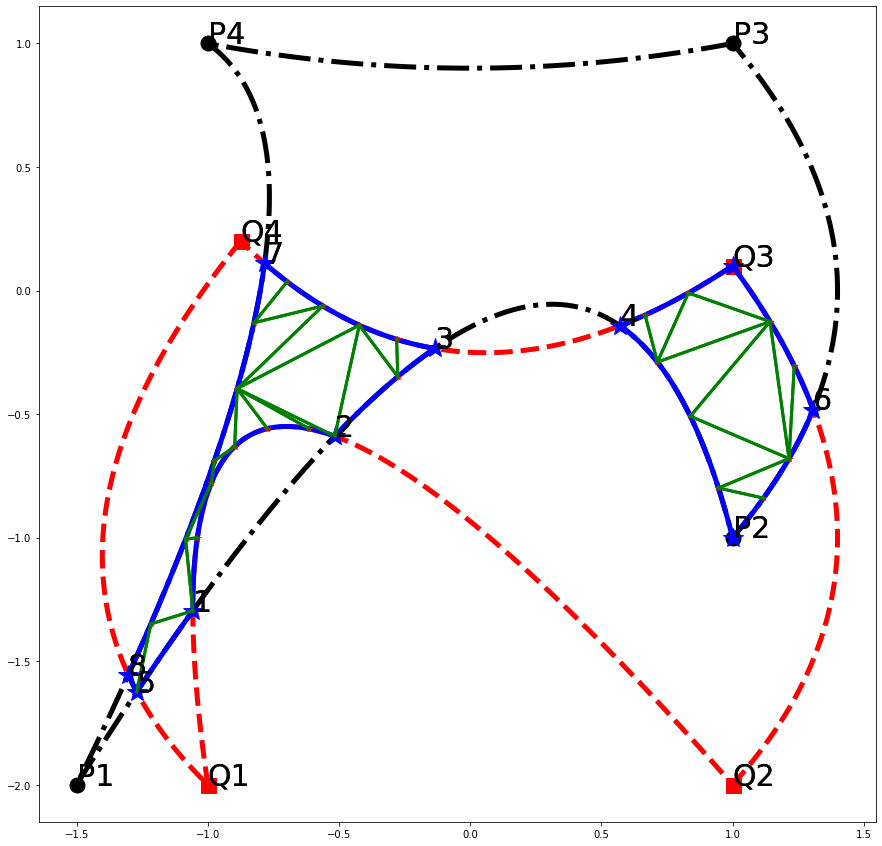}
	\includegraphics[width=7.5cm]{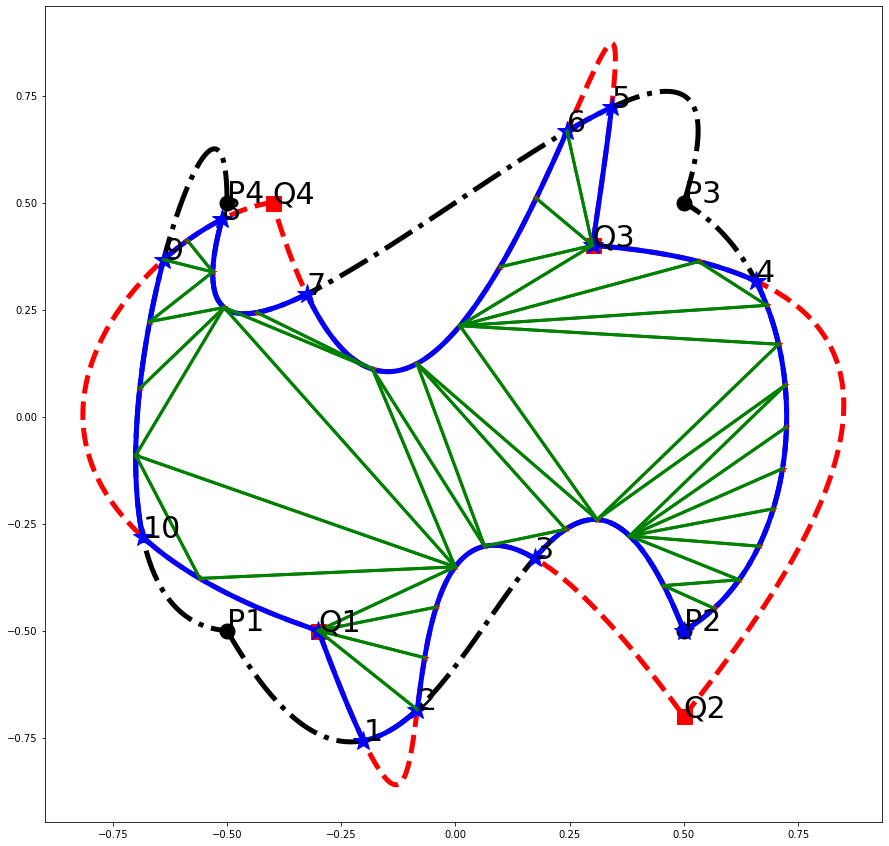}
	\caption{Partition of the intersections.}
	\label{Fig_Patition}
\end{figure}

For the two curved-edge polygons provided in Table \ref{QuadExample}, 
we compute the area of the intersections with two different approaches.
In Table \ref{Tab_Area}, the ``exact'' row denotes the exact area provided in \cite{Shashkovy2024}, 
while ``Approach A'' and ``Approach B'' represent the results obtained 
by our edge integral method \eqref{Integral_Green2} 
and triangulation-based integral method \eqref{Integral_par}, respectively. 
It can be observed that both of our methods exhibit extremely high precision, 
with integration errors close to machine zero.
Figure \ref{Fig_Patition} shows the partition of the intersection.
For the intersection enclosed by third-order isoparametric cells in Figure \ref{WA_test_3rd}, 
the results calculated using ``Approach A'' and ``Approach B'' are 1.0341187055615808
and 1.0341187055615777, respectively, which are equivalent in the sense of machine error.

\begin{table}[ht]
\centering
\caption{Area for the intersection example in Table \ref{QuadExample}. }
\begin{tabular}{c|c|c}
\toprule
\quad  & \text{Area} & Absolute error \\ \hline
\text{exact} 		& 0.723 453 730 359 014 3 	& \quad \\ \hline
\text{Approach A} 	& 0.723 453 730 359 005 6 & $8.66 \times 10^{-15} $ \\ \hline
\text{Approach B} 	& 0.723 453 730 359 006 7 & $7.55 \times 10^{-15} $ \\
\bottomrule
\end{tabular}
\label{Tab_Area}
\end{table}

\begin{remark}
In general, using Green's theorem to compute integrals over curved-edge polygons is much more efficient. 
However, this method needs to place the quadrature points along the intersection's edge, 
whereas the second approach, which triangulates the intersection and applies quadrature rules, 
ensures that the quadrature points are distributed within the interior. 
This property facilitates subsequent positivity-preserving modification. 
\end{remark}

\bigskip

\section{The Positivity-preserving modification}
\label{sec6}
\setcounter{table}{0}
\setcounter{figure}{0}
\setcounter{equation}{0}
When solving fluid dynamics with the ALE method, 
it is essential to preserve positivity for specific physical variables, such as density and internal energy. 
This necessitates a remapping method that preserves positivity, 
which primarily requires ensuring that the integral of the polynomial over the intersection remains positive. 

In our curved remapping method, we adopt the same approach as our previous work \cite{Lei2021, Lei2023}.
If the integral of a polynomial $u_i(x,y)$ over a given intersection is negative 
$\iint_{\hat{I}_{i, \tilde{i}}} u_i(x,y) dxdy < 0$, 
we adjust the polynomial based on its values at the quadrature points and the positive cell average $\bar{u}_i$.
Suppose $\hat{I}_{i, \tilde{i}}$ is partitioned into several curved-edge triangles 
$\hat{I}_{i, \tilde{i}} = \bigcup^T_{t = 1} \mathbb{T}_t$, 
and the quadrature points in $\mathbb{T}_t$ are $\{ F_{\mathbb{T}_t} (\xi_{ m}, \eta_{m} )\}_{m=1}^M$.
Define the group of quadrature points in $\hat{I}_{i, \tilde{i}}$ as
$$
G(\hat{I}_{i, \tilde{i}}) := \{  F_{\mathbb{T}_t} (\xi_{ m}, \eta_{m} ),\; m=1,\cdots,M,\; t = 1,\cdots, T  \},
$$
and the set of quadrature points in $I_i$ are $\bigcup_{\tilde{i} = 1}^{\tilde{N}} G(\hat{I}_{i, \tilde{i}}) $.

Next, we can apply the positivity-preserving limiter in \cite{Zhang2010},
compressing the polynomial $u_{i}(x, y)$ towards its positive cell average 
$\bar{u}_{i} \geq \varepsilon = 10^{-14}$,
\begin{equation}
\begin{aligned}
\hat{u}_{i} (x, y) &= \theta u_{i}(x, y) + (1-\theta) \bar{u}_{i}, \\
\theta &= \min\left\{ 1, \frac{ |\bar{u}_{i} - \varepsilon| }
{ |\bar{u}_{i} - m_{i} | } \right\} ,\quad 
m_{i} = \min\limits_{(x,y)\in \bigcup_{\tilde{i} = 1}^{\tilde{N}} G(\hat{I}_{i, \tilde{i}}) } u_{i}(x, y) ,
\end{aligned}
\end{equation}
such that the modified polynomial $\hat{u}_{i}(x, y) \geq \varepsilon > 0$ for all of the quadrature points in $I_i$.
In the meantime, this positivity-preserving modification is also conservative 
and maintains the original high-order accuracy as proved in \cite{Zhang2010}.

\bigskip

\section{Numerical examples}
\label{sec7}
\setcounter{table}{0}
\setcounter{figure}{0}
\setcounter{equation}{0}

\subsection{Accuracy test}
\label{sec7.1}
\begin{figure}
	\centering
	\includegraphics[width=5cm]{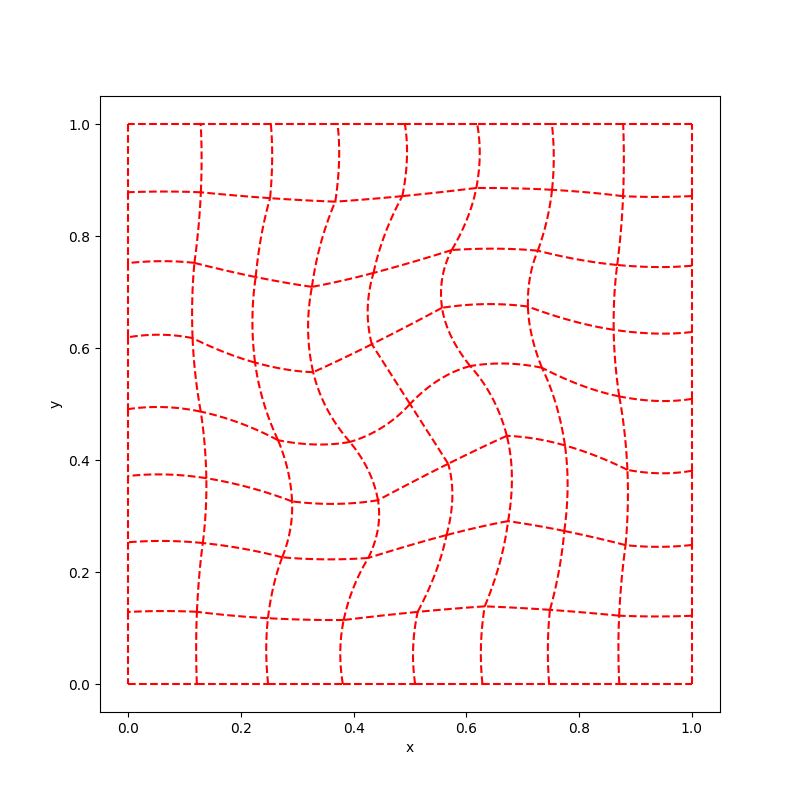}
	\includegraphics[width=5cm]{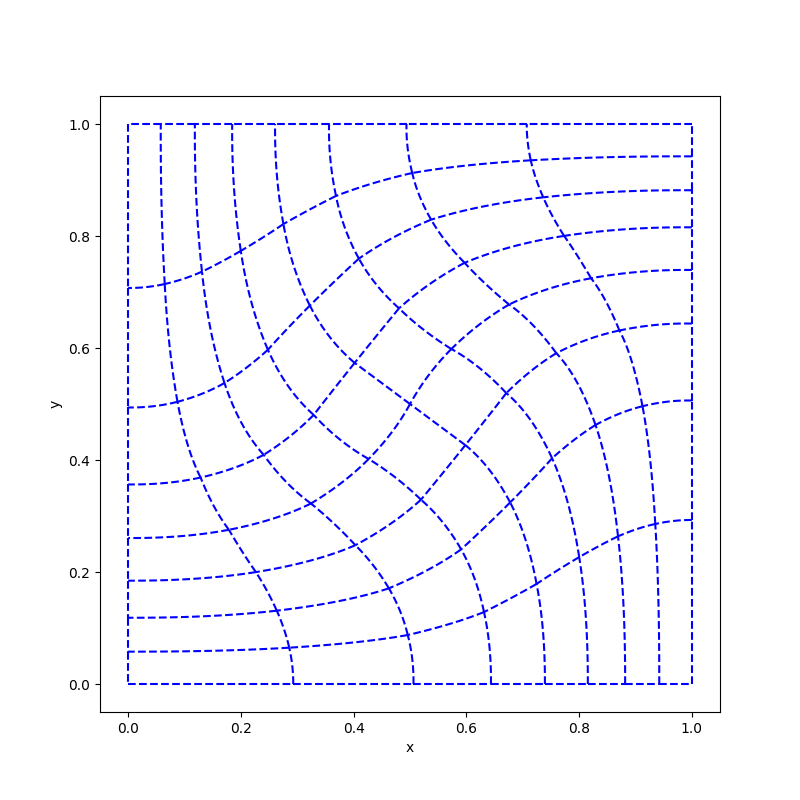}
	\includegraphics[width=5cm]{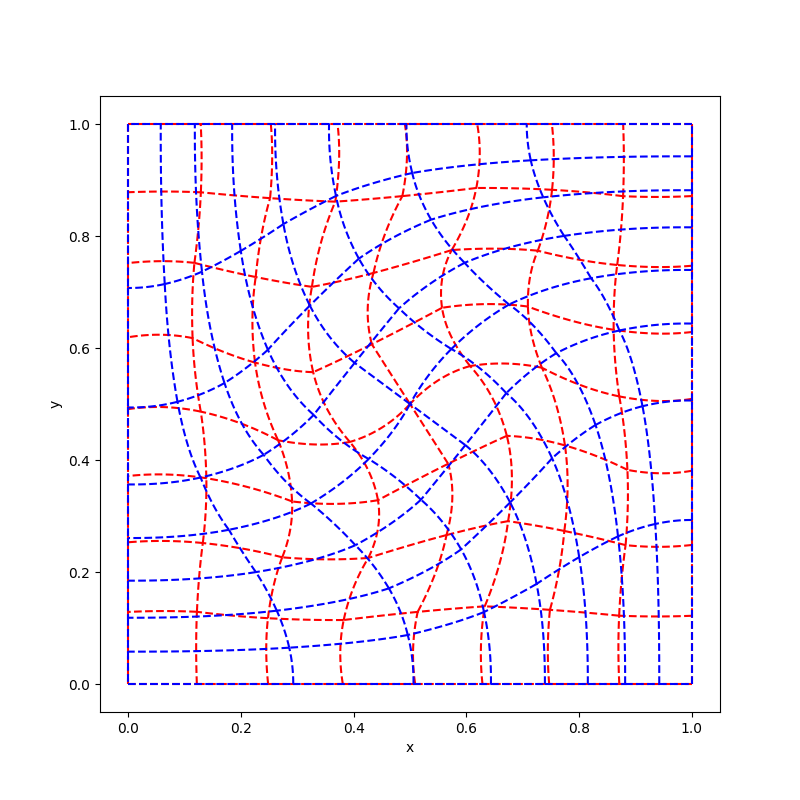}
	\caption{Left: Gresho quadratic mesh; Middle: Taylor-Green quadratic mesh; Right: meshes overlapping. }
	\label{Meshes}
\end{figure}

Following the methodology of \cite{Kenamond2021, Shashkovy2024}, 
we design a rigorous accuracy verification test for our remapping method. 
The test configuration employs two distinct curved meshes—the Gresho mesh 
and the Taylor-Green mesh—initially generated through a pure Lagrangian approach 
with straight-edge configurations for their respective benchmark problems 
(Gresho vortex and Taylor-Green flow). 
To create smooth curvilinear meshes, quadratic curves are subsequently fitted 
to all grid points derived from the Lagrangian simulation, 
as illustrated in Figure \ref{Meshes}. In this test, the Gresho mesh serves as the source domain, 
while the Taylor-Green mesh acts as the target domain, 
enabling evaluation of our curved remapping method's ability to transfer cell averages between them.

\begin{figure}
	\centering
	\includegraphics[width=7.5cm]{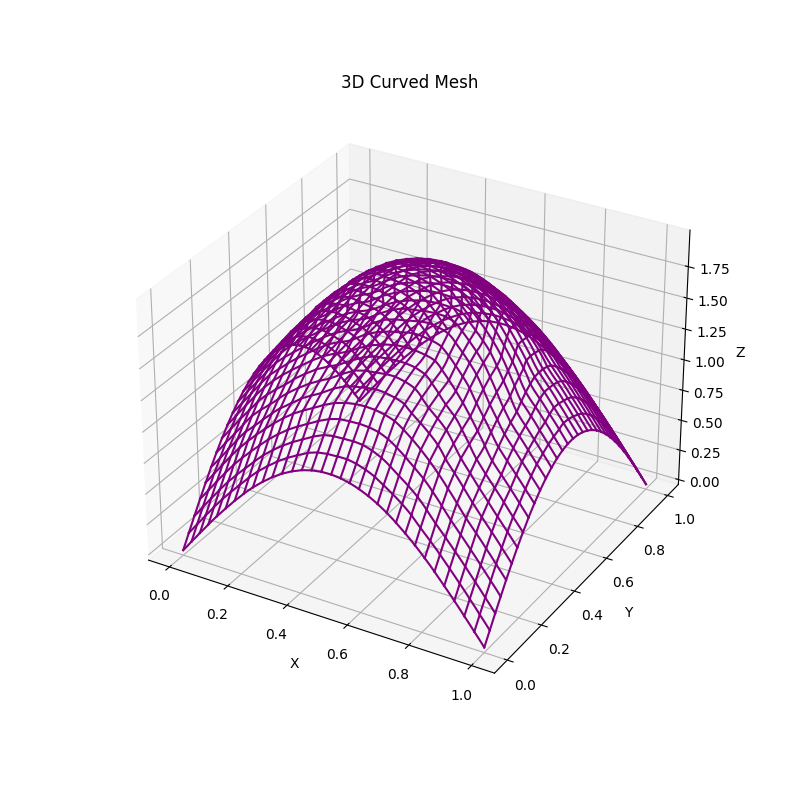}
	\includegraphics[width=7.5cm]{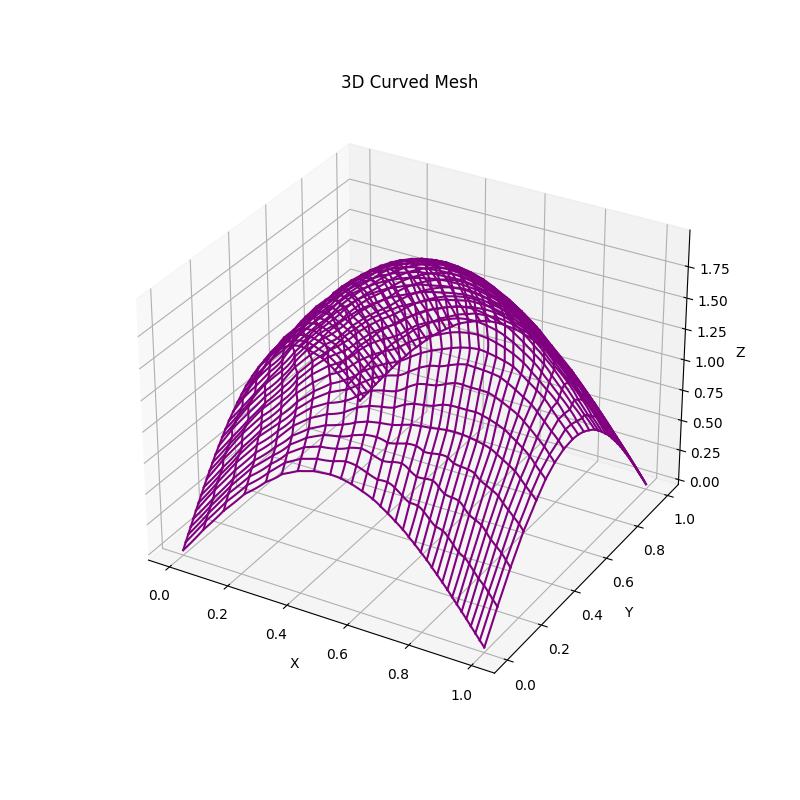}
	\caption{Left: the function \eqref{TestEqu1} on the source mesh; 
	Right: the function \eqref{TestEqu1} on the target mesh.}
	\label{Meshes3D}
\end{figure}

We define the remapped physical variable as:
\begin{equation}\label{TestEqu1}
u(x, y) = \sin( \pi x ) + \sin( \pi y ), \quad (x, y)\in [0, 1]\times [0,1].
\end{equation}
with its three-dimensional representation on 
$32\times 32$ Gresho and Taylor-Green meshes visualized in Figure \ref{Meshes3D}. 
The reference solution is established by computing exact cell averages 
directly on the target mesh, while initial values are derived from cell averages on the source mesh.

To quantify accuracy, we apply our remapping method using first-, third-, 
and fifth-order WENO reconstructions, denoted as
\( \bar{u}^{\text{1st}} \), \( \bar{u}^{\text{3rd}} \), and \( \bar{u}^{\text{5th}} \), respectively.
Error metrics are computed by comparing these reconstructed averages against  \( \bar{u}^{\text{ex}} \),
as summarized in Table \ref{Tab_Error}.  
Obviously, our remapping method maintains optimal convergence rates 
across all tested reconstruction orders, confirming its robustness in practical applications.

In Table \ref{Tab_Conservative}, we show the conservation error of our curved remapping method, 
where the clipping conservation error $e^{\text{area}}_{C}$ represents the cumulative error 
between the area of the new curved-edge cell $ |\tilde{I}_{\tilde{i}}| $ 
calculated using the WA clipping algorithm and the actual area $ |\tilde{I}^{\text{ex} }_{\tilde{i}}|$,
$$e^{\text{area}}_{C} := \sum_{\tilde{i} = 1}^{\tilde{N}} 
\left| |\tilde{I}_{\tilde{i}}| - |\tilde{I}^{\text{ex} }_{\tilde{i}}| \right|. $$
In the meantime, we define the other three conservation errors for different WENO reconstructions as 
$$
e^{\text{1st}}_{C} := \sum_{\tilde{i} = 1}^{\tilde{N}} \left| \bar{u}_{\tilde{i}}^{\text{1st}} |\tilde{I}_{\tilde{i}}| 
- \bar{u}_{\tilde{i}}^{\text{ex}} |\tilde{I}^{\text{ex} }_{\tilde{i}}| \right|, \; 
e^{\text{3rd}}_{C} := \sum_{\tilde{i} = 1}^{\tilde{N}} \left| \bar{u}_{\tilde{i}}^{\text{3rd}} |\tilde{I}_{\tilde{i}}| 
- \bar{u}_{\tilde{i}}^{\text{ex}} |\tilde{I}^{\text{ex} }_{\tilde{i}}| \right|, \;
e^{\text{5th}}_{C} := \sum_{\tilde{i} = 1}^{\tilde{N}} \left| \bar{u}_{\tilde{i}}^{\text{5th}} |\tilde{I}_{\tilde{i}}| 
- \bar{u}_{\tilde{i}}^{\text{ex}} |\tilde{I}^{\text{ex} }_{\tilde{i}}| \right|.
$$
From Table \ref{Tab_Conservative}, we can observe that the clipping error is close to machine zero, 
and the total conservation error of the three WENO reconstructions are also close to machine zero, 
which validates that our curved remapping method is conservative.

\begin{table}[htbp]
\centering
\caption{Error and order of the remapping results.}
\begin{tabular}{c|c|c|c|c|c|c}
\toprule
Mesh size  & $L^1$ error& order & $L^2$ error& order & $L^\infty$ error& order \\ \hline
\multicolumn{7}{c}{ First-order reconstruction $||\bar{u}^{\text{1st}} - \bar{u}^{\text{ex}}||$ } \\ \hline
$8\times 8$	&3.4758E-02	&		&4.3376E-02	&		&1.1057E-01	&      \\
$16\times 16$	&1.3987E-02	&1.31		&1.8662E-02	&1.22	&6.4675E-02	&0.77\\
$32\times 32$	&5.5765E-03	&1.33	&7.8870E-03	&1.24	&3.2447E-02	&1.00\\
$64\times 64$	&2.7064E-03	&1.04	&3.8512E-03	&1.03	&1.6289E-02	&0.99\\ \hline
\multicolumn{7}{c}{ Third-order WENO reconstruction $||\bar{u}^{\text{3rd}} - \bar{u}^{\text{ex}}||$ } \\ \hline
$8\times 8$	&1.0244E-03	&		&1.2736E-03	&		&3.1977E-03	&	 \\
$16\times 16$	&1.0519E-04	&3.28	&1.3841E-04	&3.20	&4.8412E-04	&2.72 \\
$32\times 32$	&1.0762E-05	&3.29	&1.4833E-05	&3.22	&6.1104E-05	&2.99 \\
$64\times 64$	&1.3221E-06	&3.03	&1.8393E-06	&3.01	&7.6675E-06	&2.99 \\ \hline
\multicolumn{7}{c}{ Fifth-order WENO reconstruction $||\bar{u}^{\text{5th}} - \bar{u}^{\text{ex}}||$ } \\ \hline
$8\times 8$	&3.0942E-05	&		&3.7514E-05	&		&9.1281E-05	&		\\
$16\times 16$	&8.3041E-07	&5.22	&1.0996E-06	&5.09	&4.0419E-06	&4.50 \\
$32\times 32$ &2.1351E-08	&5.28	&2.9594E-08	&5.22	&1.2829E-07	&4.98 \\
$64\times 64$ &6.6991E-10	&4.99	&9.4255E-10	&4.97	&4.0452E-09	&4.99 \\
\bottomrule
\end{tabular}
\label{Tab_Error}
\end{table}

\begin{table}[htbp]
\centering
\caption{Conservation error of the remapping results.}
\begin{tabular}{c|c|c|c|c}
\toprule
Mesh size  & $e^{\text{area}}_{C}$ &  $e^{\text{1st}}_{C}$ & $e^{\text{3rd}}_{C}$ & $e^{\text{5th}}_{C}$ \\ \hline
$8\times 8$	&9.9920E-16 &6.6613E-16 & 3.5527E-15	&3.7748E-15\\
$16\times 16$	&4.4409E-16 &2.2204E-16 & 4.6629E-15 & 7.3275E-15\\
$32\times 32$	&1.0103E-14 &1.0880E-14 & 2.4647E-14 & 1.1524E-13\\
$64\times 64$	&2.3204E-14 & 2.8200E-14 &5.5511E-14 & 1.5699E-13\\
\bottomrule
\end{tabular}
\label{Tab_Conservative}
\end{table}

\smallskip

\subsection{The cone test}
\label{sec7.2}
\begin{figure}
	\centering
	\includegraphics[width=7.5cm]{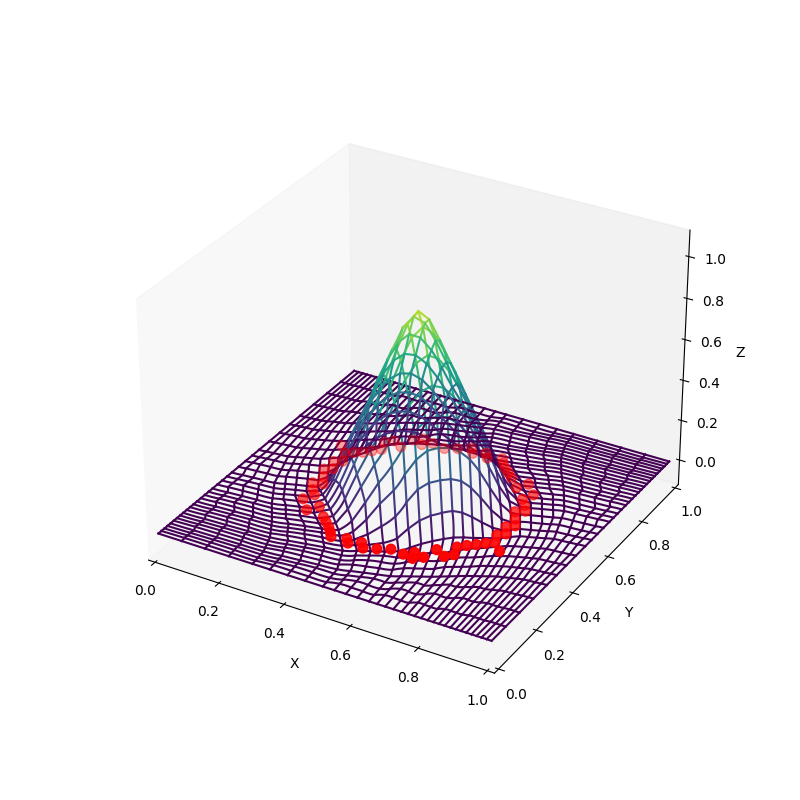}
	\includegraphics[width=7.5cm]{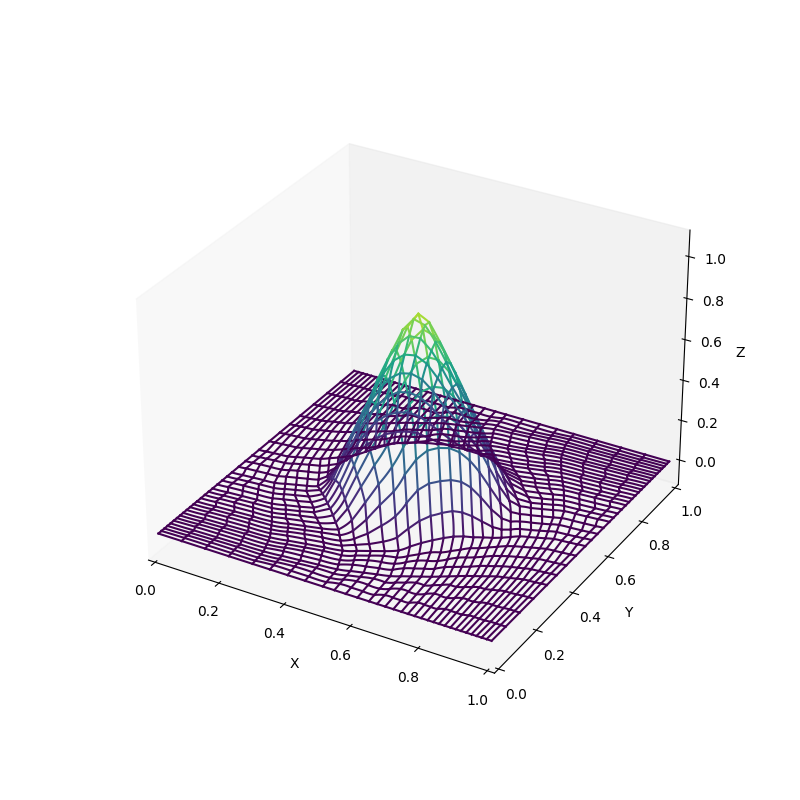}
	\caption{The cone test. Left: without the positivity-preserving limiter; 
	Right: with the positivity-preserving limiter. 
	The red markers indicate that the cell averages are negative. }
	\label{Pos_Cone3D}
\end{figure}

In this subsection, we consider the following continuous ``cone'' function, 
\begin{equation}\label{Equ_cone}
u(x, y) = \left\{ \begin{matrix} 
1 - 4r + 10^{-10}, & \sqrt{(x-0.5)^2+(y-0.5)^2} < 0.25, \\
10^{-10}, & \text{else}, \\
\end{matrix} \right.
\end{equation}
for the positivity-preserving test.
As in the previous subsection, we remap the cell averages from the source Gresho mesh
to the target Taylor-Green mesh
with the third-order WENO reconstruction.
We compare the remapping results with or without the positivity-preserving limiter, 
and we can observe from Figure \ref{Pos_Cone3D} 
that our limiter preserves positivity well.

\smallskip

\subsection{The cylinder test}
\label{sec7.3}
\begin{figure}
	\centering
	\includegraphics[width=7.5cm]{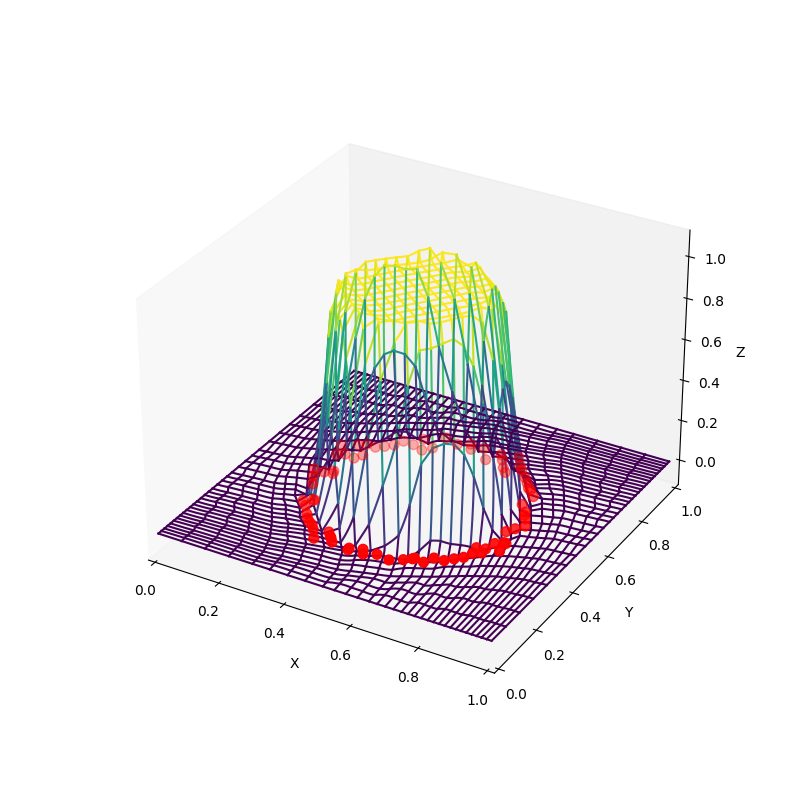}
	\includegraphics[width=7.5cm]{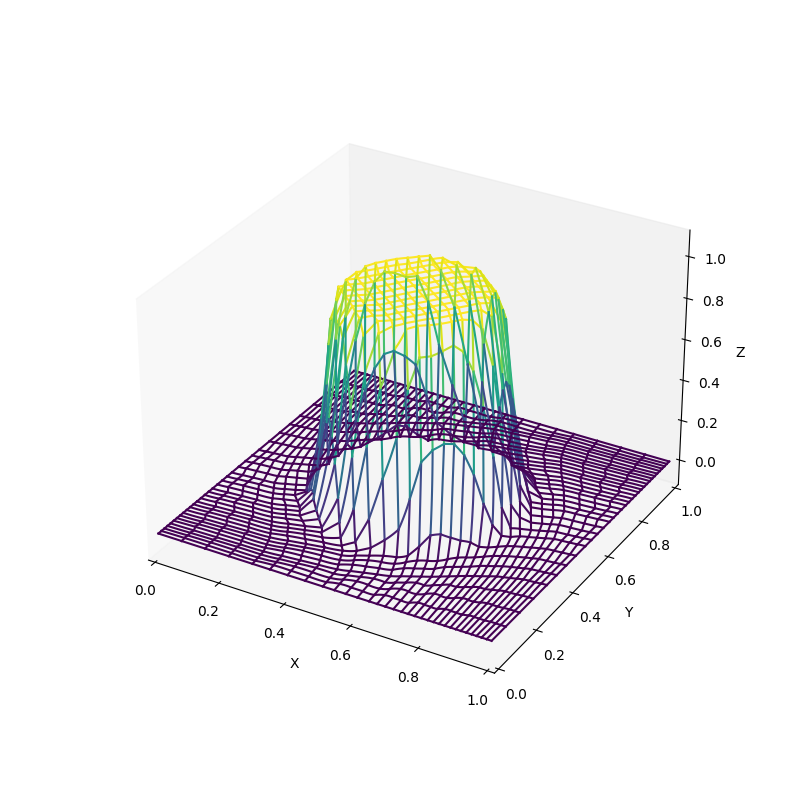}
	\caption{The cylinder test. Left: without the positivity-preserving limiter; 
	Right: with the positivity-preserving limiter. 
	The red markers indicate that the cell averages are negative. }
	\label{Pos_Cylinder3D}
\end{figure}

Next, we remap the cell averages from the source Gresho mesh
to the target Taylor-Green mesh with the third-order WENO reconstruction for
the following discontinuous ``cylinder'' function, 
\begin{equation}\label{Equ_cylinder}
u(x, y) = \left\{ \begin{matrix} 
1, & \sqrt{(x-0.5)^2+(y-0.5)^2} < 0.25, \\
10^{-10}, & \text{else}. \\
\end{matrix} \right.
\end{equation}
We conducted a comparison of the remapping results both with and without 
the positivity-preserving limiter. Upon examination of Figure \ref{Pos_Cylinder3D}, 
it is evident that our limiter effectively maintains positivity 
and that the WENO reconstruction avoids numerical oscillations in the vicinity of discontinuities.

\smallskip

\subsection{The solid body rotation test}
\label{sec7.4}
\begin{figure}
	\centering
	\includegraphics[width=7.5cm]{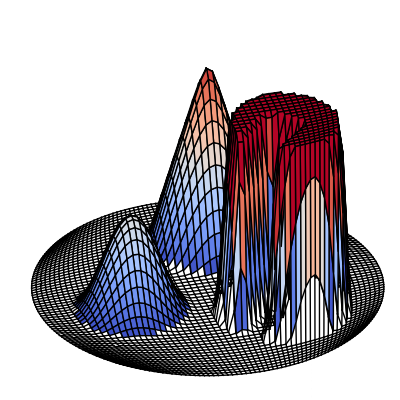}
	\includegraphics[width=7.5cm]{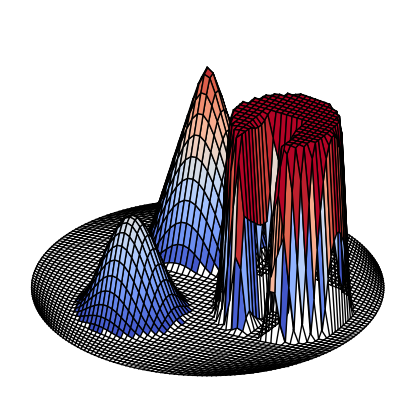}
	\includegraphics[width=7.5cm]{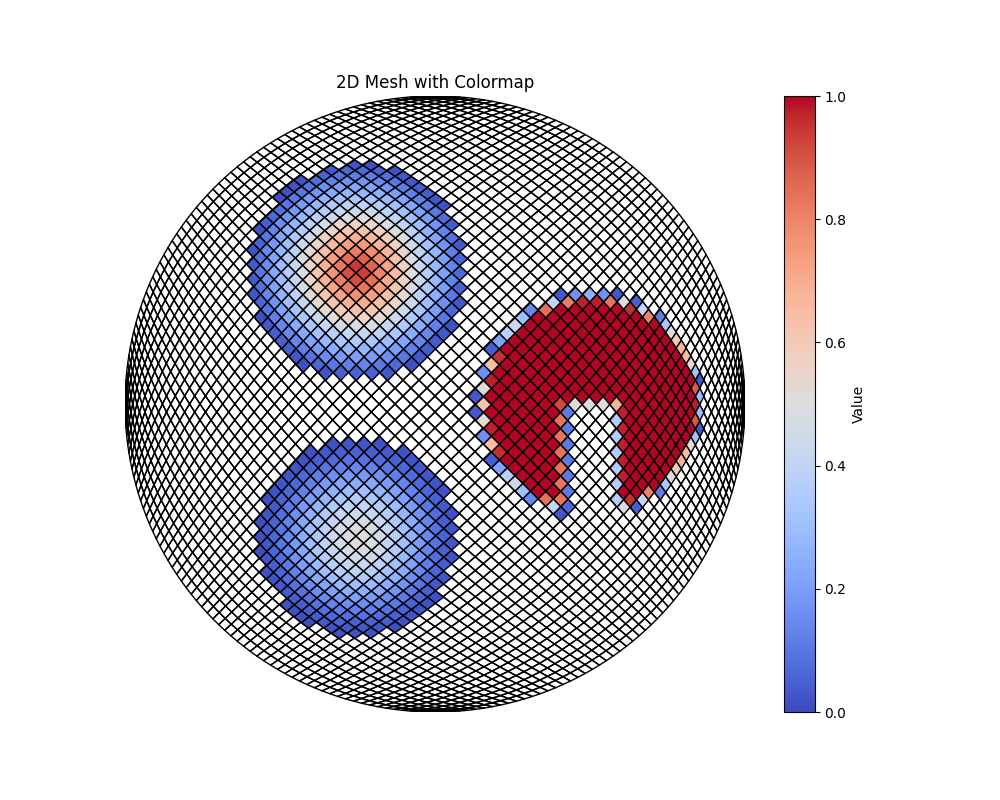}
	\includegraphics[width=7.5cm]{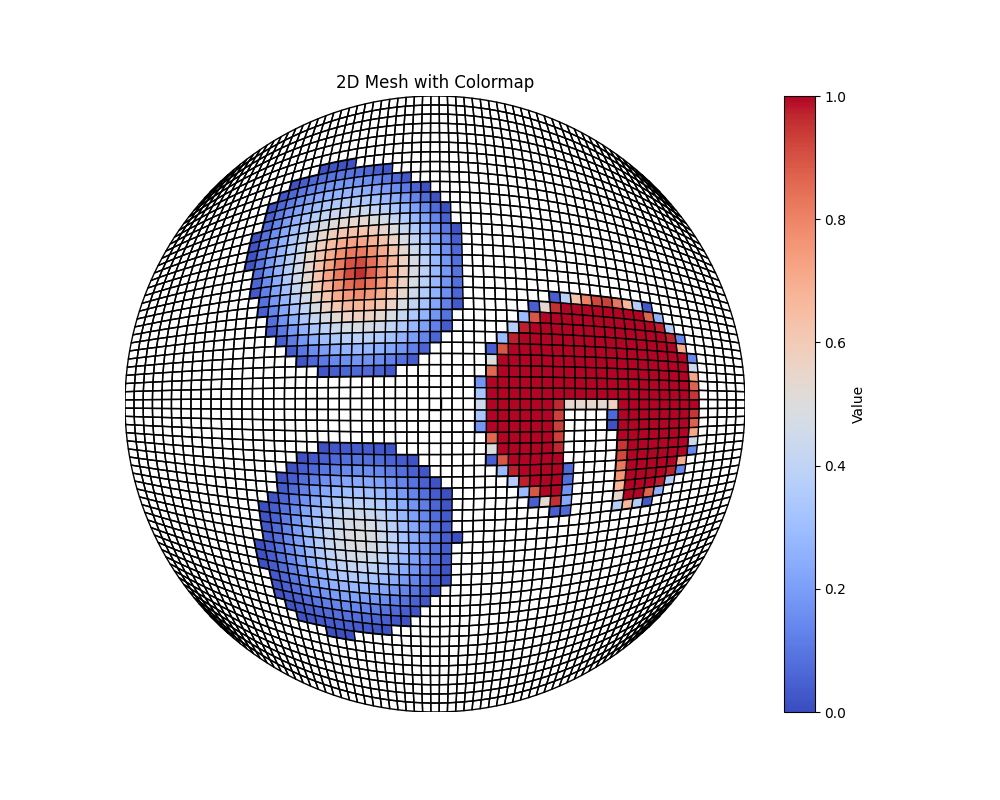}
	\caption{The solid body rotation test. Left: Initial function. Right: Rotated $\frac{\pi}{4}$ counterclockwise.}
	\label{Mixed_1}
\end{figure}

In this part, we consider the test in \cite{Anderson2015, Lipnikov2019}. 
The initial function is defined on the unit disk $x^2 + y^2 \leq 1$, 
\begin{align*}
u(x,y) = \text{cone}(x,y,\vec{r}_1) + \text{hump}(x,y,\vec{r}_1) + \text{cylinder}(x,y,\vec{r}_1) + 10^{-10}, \\
\end{align*}
where 
\begin{equation}
\begin{aligned}
\text{cone}(x, y, \vec{r}_1) & = \begin{cases}0, & \text { otherwise }  \\
1-r(x, y, \vec{r}_1), & r(x, y, \vec{r}_1) \leq 1 \end{cases} \\
\text{hump}(x, y, \vec{r}_2) & = \begin{cases}0, & \text { otherwise } \\
\frac{1}{4}(1+\cos (\pi r(x, y, \vec{r}_2)), & r(x, y, \vec{r}_2) \leq 1 \end{cases} \\
\text{cylinder}(x, y, \vec{r}_3) & = \begin{cases}0, & r(x, y, \vec{r}_3)>1 \\
0, & y < 0 \; \text{and} \; |x-\frac{1}{2}| < \frac{0.35}{4} \\
1, & \text { otherwise }\end{cases}
\end{aligned}
\end{equation}
with 
\begin{align*}
r(x, y, \vec{r}_1) &:= \frac{1}{0.35}\sqrt{(x+\frac{1}{4})^2 + (y-\frac{\sqrt{3}}{4})^2 },  \\
r(x, y, \vec{r}_2) &:= \frac{1}{0.35}\sqrt{(x+\frac{1}{4})^2 + (y+\frac{\sqrt{3}}{4})^2 },  \\
r(x, y, \vec{r}_3) &:= \frac{1}{0.35}\sqrt{(x-\frac{1}{2})^2 + y^2 }.
\end{align*}
We use Gmsh \cite{gmsh} to generate second-order curvilinear meshes with 3,481 cells.
Figure \ref{Mixed_1} displays the two-dimensional and three-dimensional visualizations of $u(x,y)$ 
evaluated on the initial mesh as well as the mesh subjected to a counterclockwise rotation of $\frac{\pi}{4}$.
Notice that, for 3D figures, we only plot the surfaces with cell averages larger than 0.01.

\begin{figure}
	\centering
	\subfigure[]{\includegraphics[width=5cm]{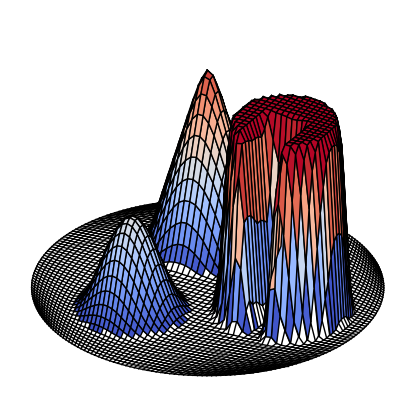}}
	\subfigure[]{\includegraphics[width=5cm]{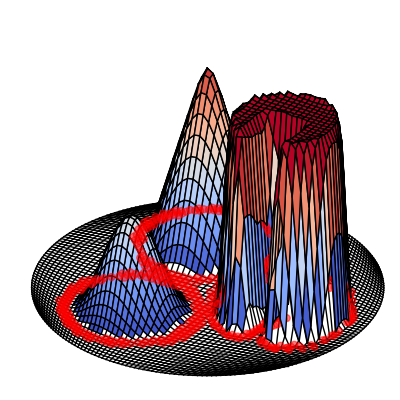}}
	\subfigure[]{\includegraphics[width=5cm]{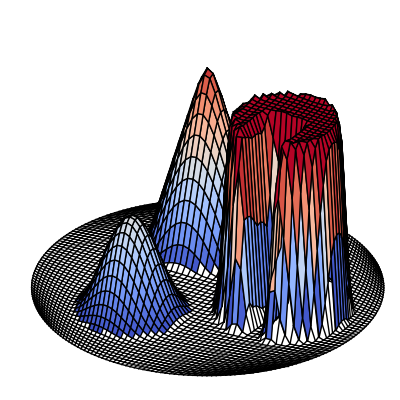}}
	\subfigure[]{\includegraphics[width=5cm]{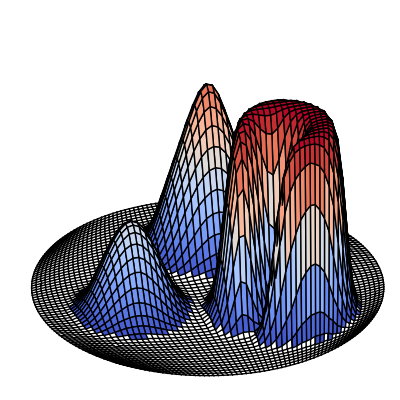}}
	\subfigure[]{\includegraphics[width=5cm]{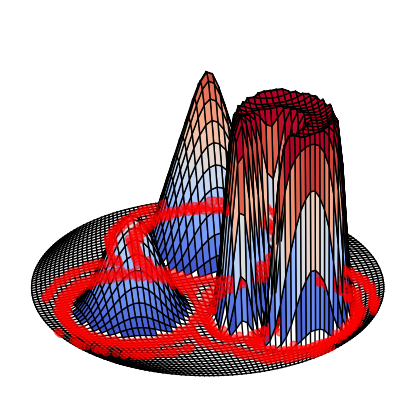}}
	\subfigure[]{\includegraphics[width=5cm]{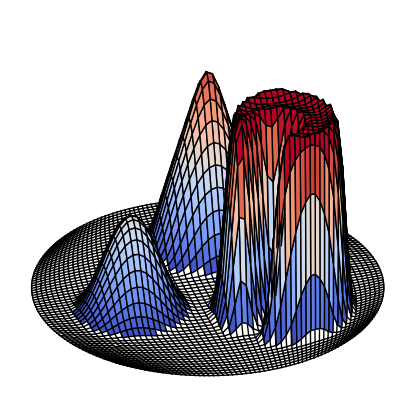}}
	\caption{The solid body rotation test.  Left column: first-order scheme.
	Middle column: third-order WENO scheme without the positivity-preserving limiter.
	Right column: third-order WENO scheme with the positivity-preserving limiter.
	First row: Rotated $\frac{\pi}{4}$ counterclockwise with remapping 1 time.
	Second row: Rotated $2\pi$ counterclockwise with remapping 8 times.	
	The red markers in the left bottom subfigure indicate that the cell averages are negative. }
	\label{Mixed3D}
\end{figure}

To test our remapping method, 
we rotate the mesh counterclockwise $2\pi$ in increments of $\frac{\pi}{4}$. 
On each segment of rotation, we apply our curved remapping method described above. 
In total, we perform 8 remapping times.
In Figure \ref{Mixed3D}, we show the remapping results
after one remapping step and eight remapping steps (returning to the initial mesh),
with first-order reconstruction, 
third-order WENO reconstruction with or without the positivity-preserving limiter.
After multiple remapping steps, our remapping method still maintains excellent results, 
demonstrating the robustness of our method. 
For the cone function, there is no smearing of the peak, and for the cylinder function, 
the high-order scheme provides sharper resolution near discontinuities without introducing significant numerical 
oscillations. Additionally, due to our intersection-based remapping method, 
we do not need to limit the rotation angle for each remapping step in this case. 
The rotation of the mesh can be achieved with very few remapping steps, 
whereas in \cite{Lipnikov2019}, 40 remapping steps were required (9 degrees per step), 
and in \cite{Anderson2015}, 1440 remapping steps were needed (1/4 degree per step).

\bigskip

\section{Conclusion}
\label{sec8}
In this paper, we present a robust intersection-based remapping method for curvilinear meshes. 
Utilizing the WA clipping algorithm, our approach efficiently handles cells of arbitrary curvature, 
even for higher-order or more complex curves, without significantly increasing computational cost. 
By incorporating high-order numerical integration techniques, 
multi-resolution WENO reconstruction, and a positivity-preserving limiter, 
our method achieves conservation, high-order accuracy, maintains essential non-oscillatory properties, 
and rigorously preserves positivity without compromising accuracy.

In the future, we aim to first extend this curved remapping method,
which is designed for finite volume schemes, to the DG framework. 
Specifically, we plan to develop high-order, conservative, 
and positivity-preserving polynomial projection curved remapping method 
tailored for DG schemes.
After that, we aim to develop a high-order curvilinear indirect ALE method, 
integrating our remapping approach to enhance accuracy and robustness in simulations 
involving large deformations, complex geometries, and multi-material interactions.

\bigskip

\normalem

\end{document}